\documentclass[12pt]{amsart}
\usepackage{amssymb}
\usepackage{amsfonts}
\usepackage{amsmath}
\usepackage{array}
\usepackage{rotating}
\usepackage{epsf}

\newtheorem{example}{Example}[section]

\newtheorem{note}[example]{Note}
\newtheorem{theorem}[example]{Theorem}

\newtheorem{corollary}[example]{Corollary}

\newtheorem{lemma}[example]{Lemma}

\def\Proof{\noindent \it Proof -- \rm}
\def\qed{\hspace{3.5mm} \hfill \vbox{\hrule height 3pt depth 2 pt width 2mm}
\bigskip}
\def\hp{.8cm}
\def\hpp{1.5cm}
\def\hdp{.4cm}

\def\FQSym{{\bf FQSym}}

\def\Std{{\rm std}}

\def\<{\langle}
\def\>{\rangle}

\def\G{{\bf G}}

\def\SG{{\mathfrak S}}

\def\Sym{{\bf Sym}}
\def\NCSF{{\bf Sym}}

\def\Des{\operatorname{Des}}

\def\PBT{{\bf PBT}}

\def\bip#1#2{\left(\begin{matrix}#1 \cr #2\end{matrix}\right)}
\def\bil#1#2{\begin{matrix}#1 \cr #2\end{matrix}}

\def\gaudend{\!\prec\!}   
\def\droitdend{\!\succ\!} 
\def\gaudendb{\ll}     
\def\droitdendb{\gg}   

\def\bas{\downarrow}
\def\droit{\triangleright}

\def\Tabvrule{\vrule width-0.4pt}       
\def\Tabhrule{\hrule \hrule height-0.4pt} 
\def\Tabstrut{\vrule height2.2ex 
                     depth0.8ex  
                     width0ex    
\relax}

\def\PasCase#1{\omit%
            $\vcenter{\hbox {\vbox to 0.4pt{}}
               \hbox{\makebox[3ex]{\Tabstrut$#1$}}}%
               \Tabvrule$}
\def\PasCasePoint{\PasCase{\cdot}}
\def\DessinCarre#1{%
    \vcenter{\hbox{}\hrule
             \hbox{\vrule\makebox[3ex]{\Tabstrut$#1$}\vrule}\Tabhrule}%
             \Tabvrule}
\def\GenRuban#1{\vcenter{\halign{&$\DessinCarre{##}$\cr#1}}\egroup}

\def\sTabvrule{\vrule width-0.4pt}
\def\sTabhrule{\hrule \hrule height-0.4pt}
\def\sTabstrut{\vrule height1.6ex depth0.6ex width0ex \relax}
\def\sDessinCarre#1{%
    \vcenter{\hbox{}\hrule
             \hbox{\vrule\makebox[2.3ex]%
                  {\sTabstrut$\scriptstyle#1$}\vrule}\sTabhrule}%
             \sTabvrule}
\def\sGenRuban#1{\vcenter{\halign{&$\sDessinCarre{##}$\cr#1}}\egroup}

\def\ruban{%
  \bgroup
  \let\ =\omit
  \let\\=\cr
  \let\x=\times
  \let\.=\PasCasePoint
  \offinterlineskip
  \GenRuban}

\def\sruban{%
  \bgroup
  \let\ =\omit
  \let\x=\times
  \let\\=\cr
  \offinterlineskip
  \sGenRuban}

\title{A one-parameter family of dendriform identities}

\author[J.-C.~Novelli and J.-Y.~Thibon]%
{Jean-Christophe Novelli and Jean-Yves Thibon}

\address[] {Institut Gaspard-Monge, Universit\'e Paris-Est\\
5 Boulevard Descartes \\Champs-sur-Marne \\77454 Marne-la-Vall\'ee cedex 2 \\
FRANCE}
\email[Jean-Christophe Novelli]{novelli@univ-mlv.fr}
\email[Jean-Yves Thibon]{jyt@univ-mlv.fr} 

\begin{document}

\begin{abstract}
We prove a $q$-identity in the dendriform dialgebra of colored free
quasi-symmetric functions. For $q=1$, we recover identities due to
Ebrahimi-Fard, Manchon, and Patras, in particular the noncommutative
Bohnenblust-Spitzer identity.
\end{abstract}

\maketitle

\section{Introduction}

The classical Spitzer and Bohnenblust--Spitzer identities~\cite{Spi,Bax,Rota}
from probability theory can be formulated in terms of certain algebraic
structures known as commutative Rota-Baxter algebras.
Recently, Ebrahimi-Fard et al.~\cite{EMP} have extended these identities to
noncommutative Rota-Baxter algebras.
Their results can in fact be formulated in terms of dendriform
dialgebras~\cite{EMP2}, a class of associative algebras whose multiplication
split into two operations satisfying certain compatibility
relations~\cite{LR1}.
Here, we exploit a natural embedding of free dendriform dialgebras into free
colored quasisymmetric functions in order to simplify the calculations, and to
obtain a $q$-analog of the main formulas of~\cite{EMP,EMP2}.

{\footnotesize
\bigskip {\it Notations -- } We  assume that the reader is familiar with
the standard notations of the theory of noncommutative symmetric functions
\cite{NCSF1,NCSF6}.
}

\section{Dendriform algebras and free quasi-symmetric functions}

\subsection{Dendriform algebras}

A dendriform dialgebra~\cite{Lod} is an associative algebra $A$ whose
multiplication $\cdot$ can be split into two operations
\begin{equation}
a\cdot b = a \gaudend b + a \droitdend b
\end{equation}
satisfying
\begin{equation}
(x\gaudend y)\gaudend z = x\gaudend (y\cdot z)\,,
(x\droitdend y)\gaudend z = x\droitdend (y\gaudend z)\,,
(x\cdot y)\droitdend z = x\droitdend (y\droitdend z)\,.
\end{equation}


\subsection{Free quasi-symmetric functions}

For example, the algebra of free quasi-symmetric functions
$\FQSym$~\cite{NCSF6} (or the Malvenuto-Reutenauer Hopf algebra of
permutations~\cite{MR}) is dendriform. Recall that for a
totally ordered alphabet $A$,  $\FQSym(A)$ is the algebra spanned
by the noncommutative polynomials
\begin{equation}
\G_\sigma(A)  := \sum_{\genfrac{}{}{0pt}{}{w\in A^n}{\Std(w)=\sigma}} w
\end{equation}
where $\sigma$ is a permutation in the symmetric group $\SG_n$ and $\Std(w)$
denotes the standardization of the word $w$.
The multiplication rule is 
\begin{equation}
\G_\alpha \G_\beta = \sum_{\gamma\in \alpha * \beta} \G_\gamma,
\end{equation}
where the convolution $\alpha*\beta$ is~\cite{MR}
\begin{equation}
\alpha * \beta =
\sum_{\genfrac{}{}{0pt}{}{\gamma=uv}{\Std(u)=\alpha ;\, \Std(v)=\beta}}
\gamma\,.
\end{equation}

The dendriform structure of $\FQSym$ is inherited from that of the free
associative algebra over $A$, which is~\cite{NT1,NT2}
\begin{eqnarray}
u\gaudend v=
\begin{cases}
uv &\mbox{if $\max(v)\le\max(u)$}\\
0 &\mbox{otherwise},\\
\end{cases}\\
u\droitdend v=
\begin{cases}
uv &\mbox{if $\max(v)\geq\max(u)$}\\
0 &\mbox{otherwise},\\
\end{cases}
\end{eqnarray}
This yields
\begin{equation}
\G_\alpha \G_\beta = \G_\alpha \gaudend \G_\beta + \G_\alpha \droitdend
\G_\beta\,,
\end{equation}
where
\begin{equation}
\G_\alpha \gaudend \G_\beta =
\sum_{\genfrac{}{}{0pt}{}{\gamma=uv \in \alpha*\beta}%
{|u|=|\alpha| ;\, \max(v)<\max(u)}}
\G_\gamma\,,
\end{equation}
\begin{equation}
\G_\alpha \droitdend \G_\beta =
\sum_{\genfrac{}{}{0pt}{}{\gamma=u.v\in \alpha*\beta}%
{|u|=|\alpha| ;\, \max(v)\geq\max(u)}} \G_\gamma\,.
\end{equation}
Then $x=\G_1$ generates a free dendriform dialgebra in $\FQSym$, isomorphic to
$\PBT$, the Loday-Ronco algebra of planar binary trees~\cite{LR1}.

There is a Hopf embedding $\iota:\ \Sym\rightarrow \PBT$ of noncommutative
symmetric functions into $\PBT$~\cite{NCSF1,NCSF6,HNT}, which is given by
\begin{equation}
\iota(S_n)=
(\dots((x\droitdend x)\droitdend x)\dots)\droitdend x\quad\text{($n$ times).}
\end{equation}
One of the identites of \cite{EMP} amounts to an expression of $\iota(\Psi_n)$
in terms of the dendriform operations. It reads
\begin{equation}
\iota(\Psi_n)=
(\dots((x\triangleright x)\triangleright x)
\dots\triangleright x \quad\text{($n$ times),}
\end{equation}
where $a\triangleright b = a\droitdend b -b\gaudend a$.
Interestingly enough, applying this identity to the dendriform products of a
Rota-Baxter algebra yields a closed form solution to the Bogoliubov recursion
in quantum field theory \cite{EMP}. 
However, using the embedding in $\FQSym$, the proof of this identity is
remarkably simple. Indeed, 
\begin{equation}
\G_\sigma\droitdend x = \G_{\sigma\cdot(n+1)}\ \text{and}\
x \gaudend \G_\sigma = \G_{(n+1)\cdot \sigma}\,,
\end{equation}
so that $\iota(S_n)=\G_{12\dots n}$.
In terms of permutations, this is therefore the standard embedding of $\NCSF$
in $\FQSym$ as the descent algebra, for which, identifying $\G_\sigma$ with
$\sigma$,
\begin{equation}
\Psi_n = [[\dots,[1,2],\dots,n-1],n].
\end{equation}
We then clearly have
\begin{align}
x\triangleright x = \G_{12}-\G_{21}=R_2-R_{11}=\Psi_2\,\\
\Psi_2\triangleright x
 = \G_{123}-\G_{213}-\G_{312}+\G_{321}=R_3-R_{21}+R_{111}=\Psi_3\,\\
\Psi_{n-1}\triangleright x
 =\G_{12\dots n}-\dots\pm \G_{n\dots 21}
 =\sum_{k=0}^{n-1}(-1)^kR_{1^k,n-k}=\Psi_n\,.
\end{align}

\subsection{Colored free quasi-symmetric functions}

Similarly, the free dendriform dialgebra of $r$ generators $x_1,\dots,x_r$ can
be realized inside the algebra $\FQSym^{(r)}$ of free quasi-symmetric
functions of level $r$~\cite{NTcoul}. It is a straightforward generalization
of $\FQSym$, built from an $r$-colored alphabet
\begin{equation}
{\bf A} := A^{(1)}\sqcup \dots \sqcup A^{(r)}
\end{equation}
where
\begin{equation}
A^{(i)} := \{a_1^{(i)}<a_2^{(i)}<\dots \}
\end{equation}
are copies of $A$.
Writing a colored word
\begin{equation}
{\bf w} = a_{i_1}^{(u_1)} \dots a_{i_n}^{(u_n)} = (w,u)
\end{equation}
where $w=a_{i_1}\dots a_{i_n}$ is the {\em underlying word} and
$u=u_1\dots u_n$ the {\em color word}, we define
\begin{equation}
\G_{\sigma,u} := \sum_{\Std(w)=\sigma} (w,u) = \G_\sigma\otimes u.
\end{equation}
Then,
\begin{equation}
\G_{\alpha,u} \G_{\beta,v}
= \sum_{\gamma\in \alpha*\beta} \G_{\gamma,u\cdot v}
= (\G_\alpha \G_\beta) \otimes uv\,,
\end{equation}
and we have again a natural dendriform structure, in which
\begin{equation}
x_1=\G_{1,1},\ \dots,\ x_r=\G_{1,r}
\end{equation}
generate a free dendriform dialgebra.

\section{The identities}

\subsection{A $q$-analog of the $\Psi_n$ with colors}

For a color word $u$ and a permutation $\sigma$ of the same size, we introduce
the \emph{biword notation}
\begin{equation}
\label{biw}
\bip{u}{\sigma} := \G_{\sigma,u}.
\end{equation}
With any color word $u=u_1\dots u_r$, we associate a nested $q$-bracketing
\begin{equation}
\Psi^u :=
\left[ \dots \left[ \bil{u_1}{1}, \bil{u_2}{2}\right]_q \dots,\bil{u_p}{p}
\right]_q
\end{equation}
where $[x,y]_q = xy-q yx$, the multiplication of nonparenthesized biletters
being ordinary concatenation,  the result being interpreted as a linear
combination of parenthesized biwords.
For example,
\begin{equation}
\Psi^{312} = \left[ \left[ \bil{3}{1}, \bil{1}{2}\right]_q,
                    \bil{2}{3}\right]_q
           = \bip{312}{123} - q\bip{132}{213}
           - q\bip{231}{312} + q^2\bip{213}{321}.
\end{equation}
%
This is an element of the free dendriform dialgebra generated by
$x_1,\dots,x_r$:
\begin{equation}
\Psi^u = (\dots (x_{u_1}\triangleright_q x_{u_2})
          \triangleright_q x_{u_3}\dots) \triangleright_q x_{u_p}
\end{equation}
where $x\triangleright_q y = x\droitdend y - q y\gaudend x$.
For $q=1$, $\triangleright_q$ is one of the two pre-Lie products always
defined on a dendriform dialgebra.

A word $u=u_1\dots u_p$ is called \emph{initially dominated} if $u_1>u_i$ for
all $i>1$. Each word has a unique increasing factorization into initially
dominated words $u^{(i)}$, \emph{i.e.},
\begin{equation}
u = u^{(1)} \cdot u^{(2)} \cdots u^{(p)}
\end{equation}
such that $u^{(1)}_1 \leq u^{(2)}_1\leq \dots \leq u^{(p)}_1$.

With a permutation $\sigma\in\SG_n$ regarded as a word with increasing
factorization $\sigma = u^{(1)}\cdots u^{(p)}$, we associate the following
element of $\FQSym^{(n)}$:
\begin{equation}
{\mathbf \Psi}^\sigma = \Psi^{u^{(1)}} \Psi^{u^{(2)}} \cdots \Psi^{u^{(p)}}.
\end{equation}
For $q=1$, this reduces to $T_\sigma(x_1,\dots,x_n)$ in the notation
of~\cite{EMP}.
Our aim is to compute the equivalent of the sum of all $T_\sigma$ in our
context.

\begin{example}
{\rm
With $n=3$, one has
}
{\scriptsize
\begin{equation}
\label{123}
{\mathbf \Psi}^{123} = \Psi^{1} \Psi^{2} \Psi^{3} =
\bip{123}{123} + \bip{123}{132} + \bip{123}{213} + \bip{123}{231} +
\bip{123}{312} + \bip{123}{321}.
\end{equation}
\begin{equation}
{\mathbf \Psi}^{132} = \Psi^{1} \Psi^{32} =
   \bip{132}{123} +  \bip{132}{213} +  \bip{132}{312}
-q \bip{123}{132} -q \bip{123}{231} -q \bip{123}{321}.
\end{equation}
\begin{equation}
{\mathbf \Psi}^{213} = \Psi^{21} \Psi^{3} =
   \bip{213}{123} +  \bip{213}{132} +  \bip{213}{231}
-q \bip{123}{213} -q \bip{123}{312} -q \bip{123}{321}.
\end{equation}
\begin{equation}
{\mathbf \Psi}^{231} =
   \Psi^{2} \Psi^{31} =
   \bip{231}{123} +  \bip{231}{213} +  \bip{231}{312}
-q \bip{213}{132} -q \bip{213}{231} -q \bip{213}{321}.
\end{equation}
\begin{equation}
{\mathbf \Psi}^{312} =
   \Psi^{312} =
   \bip{312}{123} -q \bip{132}{213} -q \bip{231}{312} +q^2 \bip{213}{321}.
\end{equation}
\begin{equation}
\label{321}
{\mathbf \Psi}^{321} =
   \Psi^{321} =
   \bip{321}{123} -q \bip{231}{213} -q \bip{132}{312} +q^2 \bip{123}{321}.
\end{equation}
}
\end{example}

Let
\begin{equation}
\Sigma_n := \sum_{\sigma\in\SG_n} {\mathbf \Psi}^{\sigma}.
\end{equation}
Summing Equations~(\ref{123}) to~(\ref{321}), one can observe that the
coefficient of each individual biword is a power of $(-q)$ multiplied by a
power of $(1-q)$.
We shall see that this is true in general.
By putting $q=1$ into $\Sigma_n$, one then recovers the result of~\cite{EMP},
namely that the sum of the $T_\sigma$ is equal to the sum of all colorings of
the identity permutation. To prove this fact, we first group permutations into
classes having the same coefficient.

\subsection{Grouping the permutations}

If the sizes of the factors of a permutation $\sigma$ into initially dominated
words are
$|u^{(1)}|=k_1$, $\dots$, $|u^{(p)}|=k_p$, we set
\begin{equation}
S(\sigma) := (k_1,\dots,k_p) = K,
\end{equation}
and call it the \emph{saillance composition} of $\sigma$.
The following tables represent the saillance compositions of all
permutations of $\SG_3$ and $\SG_4$.
\begin{equation}
\begin{array}{|c|c|c|c|}
\hline
3   & 21  & 12  & 111\\
\hline
\hline
312 & 213 & 132 & 123\\
321 &     & 231 &    \\
\hline
\end{array}
\end{equation}

\begin{equation}
\begin{array}{|c|c|c|c|c|c|c|c|}
\hline
4    & 31   & 22   & 211  & 13   & 121  & 112  & 1111\\
\hline
\hline
4123 & 3124 & 2143 & 2134 & 1423 & 1324 & 1243 & 1234 \\
4132 & 3214 & 3142 &      & 1432 & 2314 & 1342 &      \\
4213 &      & 3241 &      & 2413 &      & 2341 &      \\
4231 &      &      &      & 2431 &      &      &      \\
4312 &      &      &      & 3412 &      &      &      \\
4321 &      &      &      & 3421 &      &      &      \\
\hline
\end{array}
\end{equation}

\medskip
The saillance composition is similar to the \emph{descent composition}
$D(\sigma)=(i_1,\dots,i_s)$ whose parts are the sizes of the maximal
increasing factors of $\sigma$.
The \emph{descent set} $\Des(I)$ of a composition $I=(i_1,\ldots,i_s)$ is
\begin{equation}
\Des(I) = \{ i_1,i_1+i_2,\dots, i_1+\dots+i_{s-1}\}.
\end{equation}
If one writes $\Des(I)=\{d_1,\dots\}$, we then define $\Des(I)^-$ as the set
$\{d_1-1,d_2-1,\dots\}$.

For a color word $u$ of size $n$ and compositions $I$ and $J$ of $n$, set
\begin{equation}
R_{I}^u :=
\sum_{D(\tau)=I} \bip{u}{\tau},
\end{equation}
and
\begin{equation}
R_I^{(J)} := \sum_{S(\sigma)=J} R_I^\sigma
= \sum_{\genfrac{}{}{0pt}{}{D(\tau)=I}{S(\sigma)=J}}
\binom{\sigma}{\tau}.
\end{equation}

\begin{example}
{\rm
Regarding the biwords as bilinear operations, we have
\begin{equation}
R_{211}^{(13)} := \bip{1423+1432+2413+2431+3412+3421}{2134+3124+4123}.
\end{equation}
\begin{equation}
R_{31}^{(121)} := \bip{1324+2314}{2143+3142+3241+4132+4231}.
\end{equation}
}
\end{example}

We shall prove later that the sum $\Sigma_n$ is a linear combination of
$R_{I}^{(J)}$.
Note that the linear span of the $R_I^{(J)}$ is not a subalgebra of the
colored free quasi-symmetric functions. However, we shall
refer to it as the space of colored noncommutative symmetric functions.

\subsection{Other bases of colored noncommutative symmetric functions}

The expression of $\Sigma_n$ is simpler in a different basis.
Let us define the colored elementary basis $\Lambda_I^{(J)}$ by
\begin{equation}
\Lambda_I^{(J)} := \sum_{I'\succeq {\overline I}^\sim} R_{I'}^{(J)},
\end{equation}
where the sum runs over the compositions $I'$ finer than the conjugate
$\overline{I}^\sim$ of the mirror image $\overline I$ of $I$.
Note that this definition is independent of the color $J$.

For example,
\begin{equation}
\Lambda_{32}^{(41)} = R_{1121}^{(41)} + R_{11111}^{(41)}.
\end{equation}

\subsection{The main result}

We shall need a simple statistic on pairs of compositions. First recall
the two basic operations on compositions $I=(i_1,\dots,i_r)$
and $J=(j_1,\dots,j_s)$:
\begin{equation}
I\cdot J = (i_1,\dots,i_r,j_1,\dots,j_s),
\text{\ \ and\ \ }
I\droit J = (i_1,\dots,i_r+j_1,\dots,j_s).
\end{equation}
Now, let us define the $I$-decomposition of a composition $J$ as the unique
sequence of compositions
\begin{equation}
J\bas I= (J^{(1)},\dots, J^{(r)})
\end{equation}
such that $J^{(k)}\,\vDash i_k$ for all $k$ and
\begin{equation}
J= J^{(1)} \circ_1 J^{(2)} \circ_2 \dots \circ_{r-1} J^{(r)},
\end{equation}
where each $\circ_i$ is either $\cdot$ or $\droit$.

Let $I$ and $J$ be two compositions of $n$, and let
$(J^{(1)},\dots,J^{(r)})= J\bas I$. Write
$J^{(k)}=(j^{(k)}_1,\dots,j^{(k)}_{s_k})$.
Then the statistic $D(I,J)$ is
\begin{equation}
D(I,J)= n - \sum_{k=1}^r j^{(k)}_{s_k}.
\end{equation}

For example, with $I=(6,2,2,4,1,4)$ and $J=(3,2,4,3,2,5)$, one has
$J\bas I=((3,2,1)$, $(2)$, $(1,1)$, $(2,2)$, $(1)$, $(4))$.
Hence $D(I, J) = 19 - 1 - 2 - 1 - 2 - 1 - 4=8$.
The complete examples for sizes $3$ and $4$ are given in Section~\ref{exD}.

We then have the following simple expression:

\begin{theorem}
\label{expr-l}
\begin{equation}
\Sigma_n = \sum_{I,J} {(-1)^{l(I)-1} q^{D(I,J)} \Lambda_I^{(J)}}.
\end{equation}
\end{theorem}

The proof of this theorem relies on Theorem~\ref{expr-l-raff} and
will be given in Section~\ref{demos}.

One easily derives from this result the expansion of $\Sigma_n$ in terms of
the ribbon basis. Note that one can work in $\NCSF$ since the colors $J$
do not interfere with the change of basis between $\Lambda$ and $R$.
We have
\begin{corollary}
The sum of all ${\mathbf \Psi}^\sigma$ in $\FQSym^{(n)}$ is
\begin{equation}
\label{sommePsi}
\sum_{\sigma\in\SG_n} {\mathbf \Psi}^\sigma =
\sum_{I,J \,\vDash\, n} c_{IJ}(q)R_I^{(J)}
\end{equation}
where
\begin{equation}
c_{IJ}(q) =\left\{
\begin{array}{lc}
0 & \text{if } \Des(I)\backslash\Des(I)^- \not\subset \Des(J),\\
(-q)^{|\Des(I)\backslash\Des(J)|} (1-q)^{|\Des(I)\cap \Des(J)|} &
\text{otherwise.}
\end{array}\right.
\end{equation}
\end{corollary}

\begin{corollary}
For $q=1$, we recover the noncommutative Bohnenblust-Spitzer identity
of~\cite{EMP}:
\begin{equation}
\sum_{\sigma\in\SG_n} T_\sigma(u_1,\dots,u_n)
= \sum_{\sigma\in\SG_n} {\mathbf \Psi}^\sigma|_{q=1} =
\sum_{J\,\vDash\,n}{R_n^{(J)}} =
\sum_{\sigma\in\SG_n} \bip{12\dots n}{\sigma}.
\end{equation}
\end{corollary}

Indeed, $c_{IJ}(1)\not=0$ iff
$\Des(I)\cap \Des(J)=\emptyset$ and
$\Des(I)\backslash\Des(I)^-\subset \Des(J)$, that is,
$\Des(I)\subset\Des(I)^-$, so that $I=(n)$ (no descents). So the sum
simplifies to the sum of all colorings of the identity permutation.

\section{A refinement of Theorem~\ref{expr-l}}
\label{raff}
\label{demos}

The proof of Theorem~\ref{expr-l} relies on an induction similiar to
Newton's recursion for symmetric functions. We shall use this recursion to
state a refinement of the theorem to prove.
Let
\begin{equation}
P_I := \sum_{S(\gamma)=I} \Psi^\gamma.
\end{equation}

\begin{lemma}
Let $I=(n)$ be a one-part composition. Then
\begin{equation}
\label{pnR}
P_I = P_{(n)}
= \sum_{k=1}^{n} (-q)^{n-k} \sum_{J \,\vDash\, n-k} R_{1^{n-k},k}^{(J,k)}.
\end{equation}
\begin{equation}
\label{pnL}
P_{(n)} = \sum_{I,J} {(-1)^{n-l(I)} q^{D(I,J)} \Lambda_I^{(J)}},
\end{equation}
where the sum is taken over all compositions $J=(j_1,\ldots,j_k)$ and all
compositions $I$ belonging to the interval of the composition lattice for
the refinement order whose finest element is $(n-j_k+1,1^{j_k-1})$ and whose
coarsest element is $(n)$.
\end{lemma}

\Proof
Formula~(\ref{pnR}) is immediate by definition of $P_{(n)}$. The second
formula follows from the first one by a simple computation in the algebra of
noncommutative symmetric functions.
\qed

For an interval $[H,K]$ of the boolean lattice of compositions of $n$, let
\begin{equation}
\Lambda_I^{[H,K]} := \sum_{J\in [H,K]} \Lambda_I^{(J)}.
\end{equation}
Using this notation, Equation~(\ref{pnL}) can be rewritten as
\begin{equation}
\label{PnLb}
P_{(n)} = \sum_{I=(i_1,\ldots,i_p)\,\vDash\, n} (-1)^{n-l(I)}
      \sum_{k=n-i_1+1}^n q^{D(I,(n-k,k))}
      \Lambda_{I}^{[(n-k,k),(1^{n-k},k)]}.
\end{equation}

\begin{note}
\label{notePn}
{\rm
The characterization of the $\Lambda_I^{(J)}$ appearing in $P_n$ is simple: it
consists in the compositions $I$ and $J$ of $n$ such that the sum of the first
part of $I$ and the last part of $J$ is strictly greater than $n$.
}
\end{note}

We now need a very simple lemma on permutations.

\begin{lemma}
\label{lemPerms}
Let $\sigma\in\SG_n$.
Then $S(\sigma)=(l_1,\ldots,l_p)$ iff
\begin{equation}
S(\sigma_1,\ldots,\sigma_{l_1+\dots+l_{p-1}}) = (l_1,\ldots,l_{p-1})
\text{\ and\ }
\sigma_{l_1+\dots+l_{p-1}+1}=n.
\end{equation}
\end{lemma}

This lemma implies a recursion for $P_L$.
For any totally ordered color alphabet $C$ of size $n$, denote by
$P_I[C]$ the result of replacing each color $i$ by $c_i$ in $P_L$.
Then, by definition of $P_L$,
\begin{equation}
P_{(l_1,\dots,l_p)} = \sum_{|C'|=n-l_p,\ |C''|=l_p,\ n\in C''}
P_{(l_1,\dots,l_{p-1})}[C'] P_{(l_p)}[C''].
\end{equation}
This can be rewritten in the more suggestive form

\begin{equation}
\label{recP}
P_{(l_1,\dots,l_p)} = P_{(l_1,\dots,l_{p-1})} \droitdendb P_{(l_p)}.
\end{equation}
where the dendriform products $\gaudendb$ and $\droitdendb$ are defined in the
biword notation of~(\ref{biw}) by

\begin{equation}
\label{gauDend}
\bip{\alpha}{\beta} \gaudendb \bip{\alpha'}{\beta'} := 
\bip{\alpha\gaudend\alpha'}{\beta*\beta'}
\end{equation}

\begin{equation}
\label{droitDend}
\bip{\alpha}{\beta} \droitdendb \bip{\alpha'}{\beta'} := 
\bip{\alpha\droitdend\alpha'}{\beta*\beta'}
\end{equation}

Thanks to (\ref{PnLb}) and~(\ref{recP}), the evaluation of $P_L$ reduces to
the following:

\begin{lemma}
\label{PLElem}
Let $I$ and $J=(j_1,\dots,j_p)$ be two compositions of the same size and let
$I'$ be a composition of $n$. Then
\begin{equation}
\Lambda_{I}^{(j_1,\dots,j_p)} \droitdendb
\Lambda_{I'}^{[(n-k,k),(1^{n-k},k)]}
=
\Lambda_{I\cdot I'}^{[
   (j_1,\dots,j_{p-1}+j_p,1^{n-k},k),
   (j_1,\dots,j_{p-1},j_p,1^{n-k},k)
   ]}.
\end{equation}
\end{lemma}

\Proof
From the characterization of the right dendriform
product~(\ref{droitDend}), 
we just have to evaluate in $\FQSym$

\begin{equation}
\left(\sum_{S(\sigma)=J  ; \sigma\in\SG_n} \G_{\sigma}\right)
\droitdend
\left(\sum_{\tau_{m-k+1}=m ; \tau\in\SG_m} \G_\tau\right)
\end{equation}
that is, thanks to Lemma~\ref{lemPerms},
\begin{equation}
\sum_{\genfrac{}{}{0pt}{}{\rho=u.v; S(u)=S(\sigma)}{\rho_{n+m-k+1}=n+m}}
  \G_\rho
= \sum_{\genfrac{}{}{0pt}{}%
{S(\rho)=(j_1,\dots,j_{p-1},K,k)}{|K|=j_p+m-k; K_1\geq j_p}}
  \G_\rho.
\end{equation}
\qed


We can now state our main result:

\begin{theorem}
\label{expr-l-raff}
Let $L=(l_1,\ldots,l_p)$ be a composition of $n$. Then
\begin{equation}
P_L = \sum_{I,J} {(-1)^{n-l(I)} q^{D(I,J)} \Lambda_I^{(J)}},
\end{equation}
where the sum is taken over all pairs of compositions $(I,J)$ such that
\begin{itemize}
\item $I$ is finer than $L$,
\item For $k=1,\dots,p-1$, $\Des(J)\cap [d_k,d_{k+1}-1]\not=\emptyset$, where
$d_k=l_1+\dots+l_k$,
\item If $I\bas L=(I^{(1)},\dots,I^{(p}))$ and
         $J\bas L=(J^{(1)},\dots,J^{(p}))$, then, for all $k\in[1,p]$,
the sum of the first part of $I^{(k)}$ and the last part of $J^{(k)}$ is
strictly greater than $l_k$.
\end{itemize}
\end{theorem}

\Proof
First, Equation~(\ref{PnLb}) and Lemma~\ref{PLElem} imply that $P_L$ is a
linear combination of $\Lambda_I^{(J)}$.
It is also clear that the theorem holds if $L=(n)$. The result now follows by
induction, since it is obviously a multiplicity-free expansion thanks
to Lemma~\ref{PLElem} and since the characterization is the expected one thanks
to Note~\ref{notePn}.

The only point that remains to be proved is that the coefficient
$(-1)^{n-l(I)} q^{D(I,J)}$ is what is expected but this follows directly
from the fact that, following the notations of Lemma~\ref{PLElem},
\begin{equation}
D(I,J) + n-k = D(I\cdot I',K),
\end{equation}
for all $I'$ such that $I'_1>n-k$ and for all $K$ in the interval
\begin{equation}
[(j_1,\dots,j_{p-1},j_p+1^{n-k},k),
 (j_1,\dots,j_{p-1},j_p,1^{n-k},k)].
\end{equation}
\qed

\Proof[of Theorem~\ref{expr-l}]
Thanks to Theorem~\ref{expr-l-raff}, there only remains to prove that each
pair $(I,J)$ appears in the expansion of exactly one $P_L$.
Indeed, starting from $I$ and $J$, one glues a part of $I$ to the previous one
if there is no descent of $J$ in between those two parts. This gives a
composition $L$ such that $\Lambda_I^{(J)}$ appears in $P_L$ since it
satisfies all three conditions of Theorem~\ref{expr-l-raff}: the third
condition is the only one that remains to be checked. It is satisfied
with $L=I$ and this property remains true after each gluing,
by definition of the gluing.
Any composition strictly finer than $L$ and coarser than $I$ does not satisfy
the second condition, any other composition coarser than $I$ does not satisfy
the third condition. All other compositions do not satisfy the first
condition.
\qed

\section{Examples}
\label{exs}

\subsection{Expressions of $\Sigma_n$ in terms of $\Lambda$ and $R$}
\label{exSL}
\label{exD}

We have
\begin{equation}
\Sigma_2 = -q\Lambda_{2}^{(11)} - \Lambda_{2}^{(2)} + \Lambda_{11}^{(11)}
      + \Lambda_{11}^{(2)}.
\end{equation}
Arranging the coefficients into a matrix, whose row $I$ and column $J$
gives the value of $\Lambda_I^{(J)}$ in $\Sigma_n$, we have
\begin{equation}
M_2=
\begin{matrix}
\begin{tabular}{p{\hdp}p{\hp}p{\hp}}
& \ \ 2 & 11
\end{tabular}
\\[.1cm]
\begin{matrix} 2\\ 11\end{matrix}
\left(
\begin{matrix}
-1 & -q \\
1  & 1  \\
\end{matrix}
\right)
\end{matrix}
\end{equation}

To save space and for better readability, we shall rather give
the matrices of the exponent $D(I,J)$ itself, where $0$ is
represented by a dot.


\begin{equation}
D_2=
\begin{matrix}
\begin{tabular}{p{\hdp}p{\hdp}}
\ 2 & 11
\end{tabular}
\\[.1cm]
\begin{matrix} 2\\ 11\end{matrix}
\left(
\begin{tabular}{p{\hdp}p{\hdp}}
. & 1 \\
. & . \\
\end{tabular}
\right)
\end{matrix}
\qquad\qquad
D_3=
\begin{matrix}
\begin{tabular}{p{\hp}p{\hp}p{\hp}p{\hp}}
\ \ 3 & \ 21 & \ 12 & 111
\end{tabular}
\\[.2cm]
\begin{matrix} 3 \\ 21 \\ 12 \\ 111\end{matrix}
\left(
\begin{tabular}{p{\hp}p{\hp}p{\hp}p{\hp}}
.  & 2 & 1 & 2 \\
.  & . & 1 & 1 \\
.  & 1 & . & 1 \\
.  & . & . & . \\
\end{tabular}
\right)
\end{matrix}
\end{equation}

\begin{equation}
D_4=
\begin{matrix}
\begin{tabular}{p{\hp}p{\hp}p{\hp}p{\hp}p{\hp}p{\hp}p{\hp}p{\hp}}
\ \ \ 4 & \ \ 31 & \ \ 22 & \ 211 & \ \ 13 & \ 121 & \ 112 & 1111
\end{tabular}
\\[.3cm]
\begin{matrix} 4 \\ 31 \\ 22 \\ 211\\ 13\\ 121\\ 112\\ 1111\end{matrix}
\left(
\begin{tabular}{p{\hp}p{\hp}p{\hp}p{\hp}p{\hp}p{\hp}p{\hp}p{\hp}}
. &  3  &  2 &  3 & 1 & 3 & 2 & 3 \\
. &  .  &  2 &  2 & 1 & 1 & 2 & 2 \\
. &  1  &  . &  1 & 1 & 2 & 1 & 2 \\
. &  .  &  . &  . & 1 & 1 & 1 & 1 \\
. &  2  &  1 &  2 & . & 2 & 1 & 2 \\
. &  .  &  1 &  1 & . & . & 1 & 1 \\
. &  1  &  . &  1 & . & 1 & . & 1 \\
. &  .  &  . &  . & . & . & . & . \\
\end{tabular}
\right)
\end{matrix}
\end{equation}
Note that all columns of $M_n$ become equal when $q=1$.


Here are now the matrices of $\Sigma_2$, $\Sigma_3$, and $\Sigma_4$ in the
ribbon basis $R_{I}^{(J)}$.

\begin{equation}
M'_2=
\begin{matrix}
\begin{tabular}{p{\hp}p{\hp}}
\ \ 2 & \ 11
\end{tabular}
\\[.1cm]
\begin{matrix} 2\\ 11\end{matrix}
\left(
\begin{tabular}{p{\hp}p{\hp}}
$1$  & $1$   \\
.  & $1-q$ \\
\end{tabular}
\right)
\end{matrix}
\end{equation}

\begin{equation}
M'_3= \ \ \
\begin{matrix}
\begin{tabular}{p{\hpp}p{\hpp}p{\hpp}p{\hpp}}
\ \ \ \ 3 & \ 21 & \ 12 & 111
\end{tabular}
\\[.3cm]
\begin{matrix} 3 \\ 21 \\ 12 \\ 111\end{matrix}
\left(
\begin{tabular}{p{\hp}p{\hpp}p{\hpp}p{\hpp}}
$1$ & $1$       & $1$   & $1$       \\
.   & $1-q$     & .     & $1\!-\!q$     \\
.   & .         & $1\!-\!q$ & $1\!-\!q$     \\
.   & $-q(1\!-\!q)$ & .     & $(1\!-\!q)^2$
\end{tabular}
\right)
\end{matrix}
\end{equation}


\begin{equation}
M'_4=
\begin{matrix}
\begin{tabular}{p{\hdp}p{\hpp}p{\hpp}p{\hpp}p{\hpp}p{\hpp}p{\hpp}p{\hpp}}
\ 4 & 31 & 22 & 211 & 13 & 121 & 112 & \!\!\!\! 1111
\end{tabular}
\\[.3cm]
\left(
%
\begin{tabular}{p{\hdp}p{\hpp}p{\hpp}p{\hpp}p{\hpp}p{\hpp}p{\hpp}p{\hpp}}
$1$ & $1$ & $1$ & $1$ & $1$ & $1$ & $1$ & $1$ \\
. & $1\!-\!q$ & . & $1\!-\!q$ & . & $1\!-\!q$ & . & $1\!-\!q$ \\
. & . & $1\!-\!q$ & $1\!-\!q$ & . & . & $1\!-\!q$ & $1\!-\!q$ \\
.   &  $-q(1\!-\!q)$  &  . &  . & 1 & 1 & 1 & 1 \\
.   &  .          &  1 &  2 & . & 2 & 1 & 2 \\
.   &  .          &  1 &  1 & . & . & 1 & 1 \\
.   &  .          &  . &  1 & . & 1 & . & 1 \\
.   &  $q^2(1\!-\!q)$ &  . &  . & . & . & . & . \\
\end{tabular}
\right)
\end{matrix}
\end{equation}

\subsection{Expressions of $P_L$ in terms of $\Lambda$}
\label{plL}

The entry $(I,J)$ in the following matrices is the composition $L$ such that
$\Lambda_{I}^{(J)}$ appears in $P_L$.

\begin{equation}
N_2=
\begin{matrix}
\begin{tabular}{p{\hp}p{\hdp}}
\ \ \ 2 & 11
\end{tabular}
\\[.2cm]
\begin{matrix} 2\\ 11\end{matrix}
\left(
\begin{tabular}{p{\hdp}p{\hdp}}
2  & 2  \\
2  & 11 \\
\end{tabular}
\right)
\end{matrix}
\qquad\qquad
N_3=
\begin{matrix}
\begin{tabular}{p{\hpp}p{\hp}p{\hp}p{\hp}}
\ \ \ \ \ 3 & 21 & 12 & 111
\end{tabular}
\\[.3cm]
\begin{matrix} 3 \\ 21 \\ 12 \\ 111\end{matrix}
\left(
\begin{tabular}{p{\hp}p{\hp}p{\hp}p{\hp}}
3  & 3  & 3  & 3   \\
3  & 21 & 3  & 21  \\
3  & 12 & 12 & 12  \\
3  & 21 & 12 & 111 \\
\end{tabular}
\right)
\end{matrix} 
\end{equation}
\bigskip

\begin{equation}
N_4=
\begin{matrix}
\begin{tabular}{p{\hpp}p{\hp}p{\hp}p{\hp}p{\hp}p{\hp}p{\hp}p{\hp}}
\ \ \ \ \ \  4 &\ 31 & \ 22 & 211 & 13 & 121 & 112 & 1111
\end{tabular}
\\[.3cm]
\begin{matrix} 4 \\ 31 \\ 22 \\ 211\\ 13\\ 121\\ 112\\ 1111\end{matrix}
\left(
\begin{tabular}{p{\hp}p{\hp}p{\hp}p{\hp}p{\hp}p{\hp}p{\hp}p{\hp}}
4   & 4   & 4   & 4   & 4   & 4   & 4   & 4   \\
4   & 31  & 4   & 31  & 4   & 31  & 4   & 31  \\
4   & 22  & 22  & 22  & 4   & 22  & 22  & 22  \\
4   & 31  & 22  & 211 & 4   & 31  & 22  & 211 \\
4   & 13  & 13  & 13  & 13  & 13  & 13  & 13  \\
4   & 31  & 13  & 121 & 13  & 121 & 13  & 121 \\
4   & 22  & 22  & 22  & 13  & 112 & 112 & 112 \\
4   & 31  & 22  & 211 & 13  & 121 & 112 & 1111\\
\end{tabular}
\right)
\end{matrix}
\end{equation}

\medskip
{\footnotesize
{\it Acknowledgements.-}
This project has been partially supported by the grant ANR-06-BLAN-0380.
The authors would also like to thank the contributors of the MuPAD project,
and especially of the combinat part, for providing the development environment
for their research (see~\cite{HTm} for an introduction to MuPAD-Combinat).
}


\end{document}